\documentclass{article}

\usepackage[english]{babel}
\usepackage{amssymb}
\usepackage{amsthm}

\newtheorem{theorem}{Theorem}[section]
\newtheorem{definition}[theorem]{Definition}
\newtheorem{lemma}[theorem]{Lemma}
\newtheorem{proposition}[theorem]{Proposition}
\newtheorem{corollary}[theorem]{Corollary}
\newtheorem{remark}[theorem]{Remark}

\newcommand{\hh}{{\mathbb{H}}}
\newcommand{\rh}{\widehat{\mathbb{H}}}
\newcommand{\rr}{{\mathbb{R}}}
\newcommand{\nn}{{\mathbb{N}}}
\newcommand{\s}{{\mathbb{S}}}
\newcommand{\z}{{\mathcal{Z}}}
\newcommand{\f}{{\mathcal{F}}}
\newcommand{\m}{{\mathbb{M}}}
\newcommand{\mr}{{\mathcal{M}}}
\newcommand{\lft}{{\mathbb{G}}}
\newcommand{\rft}{{\mathcal{G}}}
\newcommand{\aff}{{\mathcal{A}}}
\newcommand{\I}{{\mathbb{I}}}
\newcommand{\B}{{\mathbb{B}}}

\title{\bf Regular Moebius transformations of the space of quaternions}
\author{Caterina Stoppato\footnote{Partially supported by GNSAGA of the INdAM, by PRIN ``Propriet\`a geometriche delle variet\`a reali e complesse'' and by PRIN ``Geometria Differenziale e Analisi Globale'' of the MIUR.} \\ 
\normalsize Dipartimento di Matematica ``U. Dini'', Universit\`a di Firenze \\ 
\normalsize Viale Morgagni 67/A, 50134 Firenze, Italy\\
\normalsize  stoppato@math.unifi.it\\}

\date{  }

\begin{document}
\maketitle


\begin{abstract}
Let $\hh$ be the real algebra of quaternions. The notion of regular function of a quaternionic variable recently presented by G. Gentili and D. C. Struppa developed into a quite rich theory. Several properties of regular quaternionic functions are analogous to those of holomorphic functions of one complex variable, although the diversity of the quaternionic setting introduces new phenomena. This paper studies regular quaternionic transformations.
We first find a quaternionic analog to the Casorati-Weierstrass theorem and prove that all regular injective functions from $\hh$ to itself are affine. In particular, the group $Aut(\hh)$ of biregular functions on $\hh$ coincides with the group of regular affine transformations. Inspired by the classical quaternionic linear fractional transformations, we define the regular fractional transformations. We then show that each regular injective function from $\rh = \hh \cup \{\infty\}$ to itself is a regular fractional transformation. Finally, we study regular Moebius transformations, which map the unit ball $\B = \{q \in \hh : |q| < 1 \}$ onto itself. All regular bijections from $\B$ to itself prove to be regular Moebius transformations.
\end{abstract}


\section{Introduction}\label{introduction}

Denote by $\hh$ the real algebra of quaternions, obtained by endowing $\rr^4$ with the following multiplication operation: if $1,i,j,k$ denotes the standard basis, define $$i^2 = j^2 = k^2 = -1,$$  $$ ij = -ji = k, jk = -kj = i, ki = -ik = j,$$ let $1$ be the neutral element and extend the operation by distributivity to all quaternions $q = x_0 + x_1 i + x_2 j + x_3 k$.  Over the last century, there have been several attempts to identify a class of quaternionic functions serving as the holomorphic functions do in the complex case. The best-known is due to R. Fueter, whose papers \cite{fueter1, fueter2, fueter3} delivered a very rich theory (see \cite{sudbery} for an excellent survey). Recent work in the field includes \cite{librodaniele} and references therein. A new theory of quaternionic functions has been proposed in \cite{cras, advances} and further developed in the subsequent papers \cite{open, zeros, milan, rigidity, poli} (for an overview, see \cite{survey}). This study is also the basis for a new functional calculus in a non commutative setting (see \cite{electronic, overview, jfa}).

The theory presented in \cite{cras, advances} is based on the following definition of regularity for quaternionic functions (inspired by C. G. Cullen \cite{cullen}). Consider the set $\s = \{q \in \hh : q^2 =-1\}$ of quaternionic imaginary units, which coincides with the (two-dimensional) unit sphere in the (three-dimensional) space of purely imaginary quaternions. For all imaginary units $I \in \s$, let $L_I = \rr + I \rr \simeq \mathbb{C}$ be the complex line through $0, 1$ and $I$. For all $f : \Omega \to \hh$ with $\Omega \subseteq \hh$ and for all $I \in \s$, denote $\Omega_I = \Omega \cap L_I$ and $f_I = f_{|_{\Omega_I}}$.

\begin{definition}\label{definition}
Let $\Omega$ be a domain in $\hh$. A real differentiable function $f : \Omega \to \hh$ is said to be \textnormal{regular} if, for all $I \in \s$, the restriction $f_I$ is holomorphic on $\Omega_I$, i.e. the function $\bar \partial_I f : \Omega_I \to \hh$ defined by
\begin{equation}
\bar \partial_I f (x+Iy) = \frac{1}{2} \left( \frac{\partial}{\partial x}+I\frac{\partial}{\partial y} \right) f_I (x+Iy)
\end{equation}
vanishes identically.
\end{definition}

As explained in \cite{advances}, a quaternionic power series $\sum_{n \in \nn} q^n a_n$ with $a_n \in \hh$ defines a regular function in its domain of convergence, which proves to be a ball $B(0,R) = \{q \in \hh : |q| <R\}$. In the same paper, it is proven that

\begin{theorem} 
If $f : B = B(0,R) \to \hh$ is regular then there exist quaternions $a_n \in \hh$ such that $f(q)=\sum_{n \in \nn} q^n a_n$ for all $q \in B$. In particular, $f \in C^{\infty}(B)$.
\end{theorem}

We may thus identify the set $\mathcal{D}_R$ of regular functions on a ball $B(0,R)$ with the set of quaternionic power series converging in the same ball. In \cite{advances} many basic results in complex analysis are extended to functions of this type: the identity principle, the maximum modulus principle, the Cauchy representation formula, the Liouville theorem, the Morera theorem and the Schwarz lemma. A very peculiar property of these functions is the distribution of their values on the 2-spheres of the type $x+y\s =  \{x+yI : I \in \s\}$ with $x,y \in \rr$ and $y \neq 0$ (see \cite{open,advances}).

\begin{theorem}
Let $f : B(0,R) \to \hh$ be a regular function and choose $x,y \in \rr$ such that $y \neq 0, x+y\s \subset B(0,R)$. There exist constants $b,c \in \hh$ such that $f(x+yI) = b + I c$ for all $I \in \s$. In other words, the map $I \mapsto f(x+yI)$ is affine (or constant, in which case we say that $x+y\s$ is a \emph{degenerate sphere} for $f$). In particular,  if $f$ has more than one zero in $x+y\s$ then it vanishes identically on $x+y\s$.
\end{theorem}

Let us list two other properties of regular functions, presented in \cite{open}, which will prove useful in the sequel.

\begin{theorem}[Minimum Modulus Principle]\label{minimum}
Let $f : B=B(0,R) \to \hh$ be a regular function. If $|f|$ has a local minimum point $p\in B$ then either $f(p)=0$ or $f$ is constant.
\end{theorem}

\begin{theorem}[Open Mapping Theorem]\label{open}
Let $f : B = B(0,R) \to \hh$ be a non-constant regular function and let $D_{f}$ be the union of the degenerate spheres of $f$ (called the \emph{degenerate set} of $f$). If $f$ is not constant then $D_f$ has empty interior and $f: B \setminus \overline{D_{f}}\to \hh$ is open.
\end{theorem}

Let us now say a few words about the algebraic structure of regular functions. Even though quaternionic multiplication does not preserve regularity, the set $\mathcal{D}_R = \{f : B(0,R) \to \hh, f \mathrm{\ regular} \}$ proves to be an associative real algebra when endowed with an appropriate multiplication operation (which we denote by $*$ and call \emph{regular multiplication}). This allows the detailed study of the zero-sets conducted in \cite{zeros, milan}. It is also possible to consider the quotient of two regular functions $f,g \in \mathcal{D}_R$ with respect to $*$-multiplication; we denote such a quotient by $f^{-*}*g$ and call it the \emph{regular quotient} of $f$ and $g$. This algebraic tool is introduced in \cite{poli} while studying the singularities of regular quaternionic functions. In the same paper such singularities are classified as \emph{essential}, \emph{removable} or \emph{poles} and a function which does not have essential singularities is called \emph{semiregular}. In section \ref{sectionpreliminary} we survey the properties of zeros and singularities of regular functions and describe their relation to the algebraic structure. These results are employed thoroughly in the present paper, which studies the quaternionic analogues of linear fractional transformations and Moebius transformations.

We begin this study by considering, in section \ref{sectionaffine}, the group $\aff$ of \emph{regular affine transformations} of $\hh$, namely the transformations $f(q) = qa + b$ with $a,b \in \hh, a \neq 0$. We first prove the following theorem.

\begin{theorem}[Casorati-Weierstrass]
Let $f$ be a regular function on $B(0,R) \setminus \{0\}$ and suppose $0$ to be an essential singularity for $f$. For each neighborhood $U$ of $0$ in $B(0,R)$, the set $f(U \setminus \{0\})$ is dense in $\hh$.
\end{theorem}

\noindent We then derive that

\begin{theorem}
A regular function $f: \hh \to \hh$ is injective if and only if it is a regular affine transformation. In particular, the group $Aut(\hh)$ of biregular functions $\hh \to \hh$ coincides with the group $\aff$ of regular affine transformations. 
\end{theorem}

\noindent We remark the complete analogy to the corresponding property of holomorphic functions of one complex variable. On the other hand, there is no resemblance to the case of several complex variables.

In section \ref{sectionlinearfractional} we mimic the study conducted in \cite{poincare}, associating to each invertible quaternionic matrix
$A=\left( \begin{array}{cc}
a & c \\
b & d
\end{array}\right) \in GL(2, \hh)$
the \emph{linear fractional transformation} of $\rh$
$$F_A (q) =(qc+d)^{-1}(qa+b)$$
($\rh = \hh \cup \{\infty\}$ denotes the Alexandroff compactification of $\hh$). The set $\lft$ of linear fractional transformations is a group with respect to the composition operation $\circ$ and setting $\Phi(A) = F_A$ defines a surjective group antihomomorphism $\Phi : GL(2,\hh)  \to \lft$ with kernel  $Ker(\Phi)=\left \{ t \I: t \in \rr \setminus \{0\} \right \}$ (where $\I$ denotes the $2 \times 2$ identity matrix). The group of regular affine transformations is a subgroup $\aff \leq \lft$ and $\lft$ is generated by the regular affine transformations and the \emph{reciprocal function} $\rho(q) = q^{-1}$. Despite the regularity of these generators of $\lft$, not all linear fractional transformations are regular. This is because composition does not preserve regularity.

For this reason, in section \ref{sectionregularfractional} we associate to each 
$A=\left( \begin{array}{cc}
a & c \\
b & d
\end{array}\right) \in GL(2, \hh)$  the \emph{regular fractional transformation} 
$$\f_A (q) = (qc+d)^{-*}*(qa+b).$$
By the formula $(qc+d)^{-*}*(qa+b)$ we denote the aforementioned regular quotient $f^{-*}*g$ of $f(q) = qc+d$ and $g(q) = qa +b$ (regular quotients are presented in detail in section \ref{sectionpreliminary}). The composition $\circ$ is not an operation on the set $\rft$ of regular fractional transformations, but we define an action of $GL(2,\hh)$ on regular quotients so that $\rft$ is the orbit of the identity function under the same action.

In section \ref{sectioncharacterization} we show that regular fractional transformations having a pole in $\widehat{\rr} = \rr \cup \{\infty\}$ are also linear fractional transformations and in particular homeomorphisms from $\rh$ onto itself. Conversely, we prove the following result.
\begin{theorem}
Let $f : \rh \to \rh$ be continuous and injective. If $f$ is semiregular in $\hh$, then either $f$ is a regular affine transformation or it is a regular fractional transformation with pole at a real point.
\end{theorem}

In section \ref{sectionball} we restrict our attention to those regular fractional transformations which map the quaternionic unit ball $\B$ onto itself, called \emph{regular Moebius transformations}. The set $\mr = \{f \in \rft : f(\B) = \B\}$ turns out to be the orbit of the identity function under the action of the symplectic group $Sp(1,1)$. Finally, we prove that

\begin{theorem}
All regular bijective functions $f : \B \to \B$ are regular Moebius transformations.
\end{theorem}


\section{Preliminary results}\label{sectionpreliminary}

We now run through the basic properties of the zeros and the singularities of regular functions. We already mentioned that if $f$ has more than one zero in $x+y\s = \{x+yI : I \in \s\}$ then it vanishes identically on $x+y\s$. For instance, the polynomial $f(q) = q^2+1$ vanishes identically on the 2-sphere $\s$ of imaginary units. In \cite{zeros} the zero-set is further characterized as follows.

\begin{theorem}\label{structure}
Let $f$ be a regular function on an open ball $B(0,R)$. If $f$ is not identically zero then its zero-set $\z_f$ consists of isolated points or isolated 2-spheres of the form $x + y \s$, for $x,y \in \rr, y\neq 0$. 
\end{theorem}

The study of the zero-set conducted in \cite{zeros} requires the introduction of the following operation on regular functions $f : B(0,R) \to \hh$.

\begin{definition}\label{multiplication}
Let $f, g$ be regular functions on an open ball $B = B(0,R)$ and consider their power series expansions $f(q) = \sum_{n \in \nn} q^n a_n, g(q) = \sum_{n \in \nn} q^n b_n$. We define the \textnormal{regular product} of $f$ and $g$ as the regular function $f*g : B \to \hh$ defined by
\begin{equation}
f*g(q) = \sum_{n \in \nn} q^n c_n, \ c_n = \sum_{k=0}^n a_k b_{n-k}.
\end{equation}
\end{definition}

Since no confusion can arise, we also write $f(q)*g(q)$ for $f*g(q)$.

\begin{remark}
Fix $R$ with $0<R\leq + \infty$ and let $\mathcal{D}_R$ be the set of regular functions $f:B(0,R)\to \hh$. Then $(\mathcal{D}_R,+,*)$ is an associative real algebra. 
\end{remark}

As observed in \cite{zeros}, the zeros of regular functions cannot be factored with respect to the standard multiplication of $\hh$. However, a factorization property is proven in terms of $*$-multiplication. 

\begin{theorem}\label{factorization}
Let $f : B = B(0,R) \to \hh$ be a regular function and let $p \in B$. Then $f(p) = 0$ if and only if there exists another regular function $g : B \to \hh$ such that $f(q) = (q-p) * g(q)$.
\end{theorem}

Finally, the zero set of a regular product is completely characterized in terms of the zeros of the two factors by the following result, which is proven in \cite{zeros}.

\begin{theorem}\label{zerosmultiplication}
Let $f,g$ be regular functions on an open ball $ B = B(0,R)$ and let $p \in B$. If $f(p) = 0$ then $f*g(p) = 0$, otherwise $f*g(p) = f(p)g(f(p)^{-1} p f(p))$. In particular $p$ is a zero of $f*g$ if and only if $f(p) = 0$ or $g(f(p)^{-1} p f(p))=0$. As a consequence, the ring $(\mathcal{D}_R,+,*)$ of regular functions on $B(0,R)$ is a domain.
\end{theorem}

Given these peculiar properties of the zeros, one may expect regular functions to have point singularities resembling the poles of holomorphic complex functions. In \cite{poli}, it is indeed proven that a quaternionic Laurent series $f(q) = \sum_{n \in \mathbb{Z}} q^n a_n$ defines a regular function in its domain of convergence, which is a spherical shell $A(0, R_1, R_2) = \{q \in \hh :R_1 < |q| < R_2\}$. This allows the construction of functions which are regular in a punctured ball $B(0,R) \setminus \{0\}$ and have a singularity at $0$. Moreover, any regular function on a spherical shell $A(0, R_1, R_2)$ admits a Laurent series expansion centered at $0$. The latter is a special case of the following.

\begin{theorem}\label{introexpansion}
Let $f$ be a regular function on a domain $\Omega$, let $p \in \hh$ and let $L_I$ be a complex line through $p$. If $\Omega$ contains an annulus $A_I = A(p,R_1,R_2) \cap L_I$ then there exists $\{a_n\}_{n \in \mathbb{Z}} \subseteq \hh$ such that $f_I(z) = \sum_{n \in \mathbb{Z}} (z-p)^n a_n$ for all $z \in A_I$. If, moreover, $p \in \mathbb{R}$ then $f$ extends to $A(p,R_1,R_2) \cup \Omega$ and $f(q) = \sum_{n \in \mathbb{Z}} (q-p)^n a_n$ for all $q \in A(p,R_1,R_2)$.
\end{theorem}

\begin{definition}
Let $f,p,R_1,R_2$ and $\{a_n\}_{n \in \mathbb{Z}} $ be as in theorem \ref{introexpansion} and suppose $R_1 = 0$. The point $p$ is called a {\em pole} if there exists an $n \in \mathbb{N}$ such that $a_{-m} = 0$ for all $m>n$; the minimum of such $n \in \mathbb{N}$ is called the {\em order} of the pole and denoted as $ord_f(p)$. If $p$ is not a pole for $f$ then we call it an {\em essential singularity} for $f$.
\end{definition}

Note that real singularities are completely analogous to singularities of holomorphic functions of one complex variable. As for non-real singularities, theorem \ref{introexpansion} only provides information on the complex line $L_I$ through the point $p$; we apparently cannot predict the behavior of the function in a (four-dimensional) neighborhood of $p$. In order to overcome this difficulty, we need to introduce the regular quotients mentioned in the introduction. We first associate to each regular function $f$ on a ball two other functions.

\begin{definition}\label{conjugate}
Let $f(q) = \sum_{n \in \nn} q^n a_n$ be a regular function on an open ball $B = B(0,R)$. We define the \textnormal{regular conjugate} of $f$, $f^c : B \to \hh$, as $f^c(q) = \sum_{n \in \nn} q^n \bar a_n$ and the \textnormal{symmetrization} of $f$, as $f^s = f * f^c = f^c*f$. 
\end{definition}

Note that $f^s(q) = \sum_{n \in \nn} q^n r_n$ with $r_n = \sum_{k = 0}^n a_k \bar a_{n-k} \in \rr$. Moreover, the zero-sets of $f^c$ and $f^s$ are characterized in \cite{zeros} as follows.

\begin{theorem}\label{conjugatezeros}
Let $f$ be a regular function on $B = B(0,R)$. For all $x,y \in \rr$ with $x+y\s \subseteq B$, the zeros of the regular conjugate $f^c$ on $x+y\s$ are in one-to-one correspondence with those of $f$. Moreover, the symmetrization $f^s$ vanishes exactly on the sets $x+y\s$ on which $f$ has a zero.
\end{theorem}

We are now ready for the definition of regular quotient given in \cite{poli}. Recall that we denote by $\z_h = \{q \in B : h(q) = 0\}$ the zero-set of a function $h$.

\begin{definition}\label{quotient}
Let $f,g : B = B(0,R) \to \hh$ be regular functions. The \emph{(left) regular quotient} of $f$ and $g$ is the function $f^{-*} * g$ defined in $B \setminus \z_{f^s}$ by $f^{-*} * g (q) = \frac{1}{f^s(q)} f^c * g(q)$. Moreover, the \emph{regular reciprocal} of $f$ is the function $f^{-*} = f^{-*} * 1$.
\end{definition}

Since no confusion can arise, we often write $(f(q))^{-*}$ for $f^{-*}(q)$. Regular quotients prove to be regular in their domains of definition. The algebraic meaning of regular quotients is explained by the following result, proven in \cite{poli}.

\begin{proposition}
Let $0 < R \leq + \infty$ and consider domain $(\mathcal{D}_R,+,*)$ of regular functions on $B(0,R)$. If we endow the set of left regular quotients $\mathcal{L}_R = \{f^{-*}*g : f,g \in \mathcal{D}_R, f \not \equiv 0\}$ with the multiplication $*$ defined by  $(f^{-*}*g)*(h^{-*}*k) = \frac{1}{f^{s}h^{s}} f^c*g*h^c*k$ then $(\mathcal{L}_R,+,*)$ is a division algebra over $\mathbb{R}$ and it is the classical ring of quotients of $(\mathcal{D}_R,+,*)$.
\end{proposition}

For the definition of the classical ring of quotients, see \cite{rowen}. Furthermore, the following relation between the regular quotient $f^{-*} * g(q)$ and 
the quotient $f(q)^{-1} g (q)= \frac{1}{f(q)} g(g)$ is proven in \cite{open,poli}.

\begin{theorem}\label{quotients}
Let $f,g$ be regular functions on $B=B(0,R)$. If we set $T_f(q) = f^c(q)^{-1} q f^c(q)$ for all $q \in B \setminus \z_{f^s}$, then
\begin{equation}
f^{-*}*g(q)= \frac{1}{f\circ T_f(q)} g \circ T_f(q)
\end{equation}
for all $q \in B \setminus \z_{f^s}$. For all $x,y \in \rr$ with $x+y\s\subset B \setminus \z_{f^s}$, the function $T_f$ maps $x+y\s$ to itself (in particular $T_f(x) = x$ for all $x \in \rr$). Furthermore, $T_f$ is a diffeomorphism from $B \setminus \z_{f^s}$ onto itself, with inverse $T_{f^c}$.
\end{theorem}

Regular quotients allow a detailed study of the poles. By analogy with meromorphic complex functions, we give the following definition.

\begin{definition}\label{semiregular}
A function $f$ is \emph{semiregular} if it does not have essential singularities or, equivalently, if the restriction $f_I$ is meromorphic for all $I \in \mathbb{S}$.
\end{definition}

As proven in \cite{poli}, $f$ is semiregular in $B(0,R_0)$ if and only if $f_{|_{B(0,R)}}$ is a regular quotient for all $R<R_0$. This allows the definition of a multiplication operation $*$ on the set of semiregular functions on a ball and the proof of the following result. 

\begin{theorem}
Let $f$ be a semiregular function on $B = B(0,R)$, choose $p = x+yI \in B$ and let $m = ord_f(p), n = ord_f(\bar p)$. There exists a unique semiregular function $g$ on $B$ such that
\begin{equation}
f(q) = \left [(q-p)^{*m}*(q-\bar p)^{*n}\right ]^{-*}* g(q)
\end{equation}
The function $g$ is regular near the sphere $x+y\s$ through $p$ and $\bar p$. Moreover, $g(p) \neq 0 \neq g(\bar p)$ (provided $m>0$ or $n>0$).
\end{theorem}

The previous result allows the study of the distribution of the poles.

\begin{theorem}\label{distribution}
If $f$ is a semiregular function on $B = B(0,R)$ then $f$ extends to a regular function on $B$ minus a union of isolated real points $x \in \rr$ or isolated 2-spheres of the type $x+y\mathbb{S}$ with $x,y \in \mathbb{R}, y \neq 0$. All the poles on each 2-sphere $x+y\mathbb{S}$ have the same order with the possible exception of one, which must have lesser order.
\end{theorem}


\section{Regular affine transformations}\label{sectionaffine}

For all $a, b \in \hh$ with $a \neq 0$, the function $l : \hh \to \hh$ defined by $l(q) = qa+b$ is regular and bijective. We will call such an $l$ a {\em regular affine transformation} of $\hh$. 
\begin{remark}
The set $\aff$ of regular affine transformations is a group with respect to the composition operation $\circ$.
\end{remark}
Indeed, the composition of two regular affine transformations is still regular and affine. Moreover, the inverse function of $l(q) = qa+b$ is $l^{-1}(q) = q a^{-1} - b a^{-1}$. We now show that all injective regular functions $\hh \to \hh$ are regular affine transformations. In order to prove this, we need to take some preliminary steps. The first one is the analogue of the Casorati-Weierstrass theorem.

\begin{theorem}\label{casorati}
Let $f$ be a regular function on $B(0,R) \setminus \{0\}$ and suppose $0$ to be an essential singularity for $f$. For each neighborhood $U$ of $0$ in $B(0,R)$, the set $f(U \setminus \{0\})$ is dense in $\hh$.
\end{theorem}

\begin{proof}
Suppose that for some $0<r\leq R$ there existed a $p \in \hh$ and an $\varepsilon >0$ such that $f(B(0,r) \setminus \{0\}) \cap B(p, \varepsilon) = \emptyset$. Setting $g(q) = (f(q) - p)^{-*}$ would then define a regular function $g : B(0,r) \setminus \{0\} \to \hh$ with $|g| \leq 1/\varepsilon$ since (by theorem \ref{quotients})
$$\frac{1}{|g(q)|}= |f(T_{f-p}(q)) - p| \geq \inf_{w \in B(0,r) \setminus \{0\}} |f(w) - p| \geq \varepsilon$$ 
for all $q \in B(0,r) \setminus \{0\}$. The point $0$ would then be a removable singularity for $g$, i.e. $g$ would be regular in $B(0,r)$. The function $f = g^{-*} + p$ would then be semiregular in $B(0,r)$. This is impossible, since we supposed $0$ to be an essential singularity for $f$.
\end{proof}

The second step is the following property of polynomials.

\begin{lemma}\label{polynomials}
Let $f(q) = q^n a_n + ...+ q a_1 +a_0$ be polynomial of degree $n$. If $f$ is injective, then $n=1$.
\end{lemma}

\begin{proof}
Suppose by contradiction that $n>1$. Since $f(0) = a_0$ and $f$ is injective, $f$ does not equal $a_0$ at any point other than $0$. In other words, $f(q) - a_0 = q^{n} a_n + ...+ q a_1 = q (q^{n-1} a_n + ...+ q a_2 + a_1)$ only vanishes for $q=0$. Thus $q^{n-1} a_n + ...+ a_1$ can only vanish at $0$. For $n>1$, by the fundamental theorem of algebra for quaternions (see, for instance, \cite{zeros, fundamental, niven, shapiro}), $q^{n-1} a_n + ...+ q a_2 + a_1$ must have a root. Thus it vanishes at $0$ and we conclude $a_1 = 0$. Iterating this process proves $a_j = 0$ for $1 \leq j \leq n-1$, so that $f$ reduces to $f(q) = q^n a_n + a_0$. We then find a contradiction: for $n>1$ the monomial $q^n$ is clearly not injective and neither is $q^n a_n + a_0$.
\end{proof}

We can now prove the desired result.

\begin{theorem}\label{affine}
A function $f: \hh \to \hh$ is regular and injective if and only if there exist $a,b \in \hh$ with $a \neq 0$ such that $f(q) = qa+b$ for all $q \in \hh$.
\end{theorem}

\begin{proof}
We only have to prove one implication. Take an injective regular $f: \hh \to \hh$ and consider the regular function $g : \hh \setminus \{0\} \to \hh$ defined by $g(q) = f(q^{-1})$. If $0$ were an essential singularity for $g$ then, by theorem \ref{casorati}, $g(B(0,1) \setminus \{0\})$ would be dense in $\hh$. By the open mapping theorem \ref{open}, $f(B(0,1))$ is an open subset of $\hh$ (the injectivity implies $D_f = \emptyset$). We would then have
$$ \emptyset \neq g(B(0,1) \setminus \{0\}) \cap f(B(0,1)) = f(\hh \setminus B(0,1)) \cap f(B(0,1)),$$
so that $f$ would not be injective.
Thus $0$ is a pole for $g$ and $f$ must be a polynomial. Finally, since $f$ is injective it must have degree $1$ by lemma \ref{polynomials}.
\end{proof}

\begin{corollary}
The group $Aut(\hh)$ of biregular functions on $\hh$ (namely, regular functions $\hh \to \hh$ having regular inverse) coincides with the group $\aff$ of regular affine transformations. 
\end{corollary}


\section{Linear fractional transformations}\label{sectionlinearfractional}

We continue our study by considering quaternionic analogues of the linear fractional transformations. The class of quaternionic functions $q \mapsto (aq+b)(cq+d)^{-1}$ is studied in \cite{poincare}. We now recall some of the results therein, but transposing them to the case of the functions $q \mapsto (qc+d)^{-1}(qa+b)$ in order to include all regular affine transformations $q \mapsto q a +b$. Denote by $\rh = \hh \cup \{\infty\}$ the Alexandroff compactification of $\hh$, which is homeomorphic to the (four-dimensional) unit sphere of $\rr^5$. Furthermore, denote by $GL(2,\hh)$ the group of invertible $2 \times 2$ quaternionic matrices.

\begin{definition}
For any $A =
\left( \begin{array}{cc}
a & c \\
b & d
\end{array}\right)
\in GL(2,\hh),$
the continuous function $F_A : \rh \to \rh$ defined by
\begin{equation}
F_A(q) = (qc+d)^{-1}(qa+b),
\end{equation}
is called the \textnormal{linear fractional transformation} associated to $A$.
\end{definition}

In the above definition we mean $F_A(\infty) = \infty$ if $c=0$, $F_A(-dc^{-1}) = \infty$ and $F_A(\infty) = c^{-1} a$ if $c \neq 0$. The following can be proven as in \cite{poincare}.

\begin{theorem}\label{antihomomorphism}
The set $\lft = \{F_A : A \in GL(2, \hh) \}$ is a group with respect to the composition operation $\circ$ and the map $$
\begin{array}{rcl}
\Phi : GL(2, \hh)  & \to &\lft\\
A & \mapsto & F_A
\end{array}
$$ is a surjective group antihomomorphism with kernel $Ker(\Phi) = \left \{ t \I: t \in \rr \setminus \{0\} \right \}$ (where $\I$ denotes the identity matrix). Moreover, the restriction of $\Phi$ to the special linear group $SL(2, \hh)$ is still surjective and has kernel $\{ \pm  \I\}$.
\end{theorem}

For all $A \in GL(2, \hh)$, the continuous function $F_A : \rh \to \rh$ has inverse $F_A^{-1} = F_{A^{-1}}$. Hence $F_A$ is a homeomorphism from $\rh$ onto itself. We now prove the following.

\begin{proposition}\label{generators}
The group $\aff$ of regular affine transformations is a subgroup $\aff \leq \lft$. The group $\lft$ is generated by $\aff \cup \{\rho\}$ where $\rho(q) = q^{-1}$ denotes the \emph{reciprocal function}.
\end{proposition}

\begin{proof} The first statement is obviously true, let us prove the second one. If $c \neq 0$ then setting $l(q) = qc+d$ and $l'(q) = q (b-dc^{-1}a) + c^{-1}a$ yields
$$(qc+d)^{-1}(qa+b) = l' \circ \rho \circ l(q)$$
Otherwise we set $l(q) = qa+b$ and $l'(q) = q d$ and observe that 
$$d^{-1}(qa+b) = d^{-1} l(q) = \rho \circ l' \circ \rho \circ l(q).\qedhere$$  
\end{proof}

Let us remark that proposition \ref{generators} can also be proven by means of theorem \ref{antihomomorphism}, thanks to the following three facts:
\begin{enumerate}
\item $\aff$ is the isomorphic image through $\Phi$ of the group of matrices of the type 
$\left( \begin{array}{cc}
a & 0 \\
b & 1
\end{array}\right)$ with $a \neq 0$. 

\item $\rho$ is the image through $\Phi$ of the matrix $
\left( \begin{array}{cc}
0 & 1 \\
1 & 0
\end{array}\right)$.

\item $GL(2, \hh)$ is generated by the set $\left \{ \left( \begin{array}{cc}
a & 0 \\
b & 1
\end{array}\right) : a,b \in \hh, a \neq 0 \right\} \cup \left\{
\left( \begin{array}{cc}
0 & 1 \\
1 & 0
\end{array}\right)
\right\}$.
\end{enumerate}

We conclude this section noticing that, despite the regularity of the above mentioned generators of $\lft$, not all linear fractional transformations are regular. Indeed, given a regular function $f$ and an affine transformation $l \in \aff$, say $l(q) = q \alpha + \beta$, the composition $l \circ f= f \alpha+ \beta$ is still regular, but we cannot say the same for the composition $\rho \circ f = f^{-1}$ with $\rho(q) = q^{-1}$. For instance, the linear fractional transformation $q \mapsto k^{-1}q^{-1} = (qk)^{-1}$ is not regular if $k \in \hh \setminus \rr$, even though $q \mapsto qk$ is. We shall try to ``fix'' this by considering the regular reciprocal $f^{-*}$ instead of the reciprocal $f^{-1}$.


\section{Regular fractional transformations}\label{sectionregularfractional}

Recall that, as explained in section \ref{sectionpreliminary}, the ring $(\mathcal{D}_\infty,+,*)$ of regular functions $\hh \to \hh$ admits a ring of quotients $\mathcal{L}_\infty$. The quotient of $f$ and $g$ is denoted by $f^{-*}*g$ and it is semiregular in $\hh$. 

\begin{definition}\label{regularfractional}
For any $A =
\left( \begin{array}{cc}
a & c \\
b & d
\end{array}\right)
\in GL(2,\hh),$
we define the \textnormal{regular fractional transformation} associated to $A$ as
\begin{equation}
\f_A(q) = (qc+d)^{-*}*(qa+b).
\end{equation}
Moreover, we denote $\rft = \{\f_A : A \in GL(2, \hh) \}$.
\end{definition}

Let us postpone the discussion of the domain of definition of $\f_A$ and focus on the algebraic properties of $\rft$, proving that it is an orbit of the group action defined below.

\begin{theorem}\label{orbit}
Choose $R>0$ and consider the ring of quotients of regular quaternionic functions in $B(0,R)$, denoted by $\mathcal{L}_R$.
Setting 
\begin{equation}
f.A = (f c + d)^{-*}*(f a+ b)
\end{equation} 
for all $f \in \mathcal{L}_R$ and for all $A =
\left( \begin{array}{cc}
a & c \\
b & d
\end{array}\right) \in GL(2,\hh)$
defines a right action of $GL(2,\hh)$ on $\mathcal{L}_R$. The normal subgroup $N=\left \{t \I: t \in \rr \setminus \{0\} \right \}$ of $GL(2,\hh)$ is included in the stabilizer of any element of $\mathcal{L}_R$ and the action of $GL(2,\hh)/N \cong SL(2,\hh)/\{\pm\I\} \cong \lft$ is faithful.
\end{theorem}

\begin{proof}
If $A =
\left( \begin{array}{cc}
a & c \\
b & d
\end{array}\right)$ and
$B =
\left( \begin{array}{cc}
\alpha & \gamma \\
\beta & \delta
\end{array}\right),$
then
$AB = \left( \begin{array}{cc}
a\alpha + c\beta & a \gamma + c \delta \\
b\alpha + d\beta & b\gamma + d\delta
\end{array}\right)$. We compute:
$$ (f.A).B =  \left \{(f.A) \gamma + \delta\right \}^{-*}*\left \{(f.A) \alpha+ \beta\right \} =$$
$$=\left \{(fc+d)^{-*}*\left [(fa+b) \gamma + (fc+d) \delta\right]\right \}^{-*}*\left \{(fc+d)^{-*}*\left [(fa+b) \alpha+ (fc+d)\beta\right]\right \} =$$
$$= \left [(fa+b) \gamma + (fc+d) \delta\right]^{-*}*(fc+d)*(fc+d)^{-*}*\left [(fa+b) \alpha+ (fc+d)\beta\right] =$$
$$= \left [f (a\gamma + c\delta) + b\gamma + d \delta\right]^{-*}*\left [f (a \alpha + c \beta)+ b \alpha +d\beta\right] = f.AB$$
Moreover, $f.\I = (0+1)^{-*}*(f 1+0) = 1^{-1} f = f$. The last statement is proven observing that
$$(fc+d)^{-*}*(fa+b) = f \Longleftrightarrow fa+b = (fc+d)*f \Longleftrightarrow fa+b = f c *f + d *f$$
and that the latter holds for all $f \in \mathcal{L}_R$ if, and only if $b=c=0, d=a = t \in \rr$.
\end{proof}

We notice that not all stabilizers coincide with $N = \left \{t \I: t \in \rr \setminus \{0\} \right \}$. For instance, the stabilizer of the identity function is the subgroup $\left \{a \I: a \in \hh \setminus \{0\} \right \} \leq GL(2,\hh)$, since
$$(qc+d)^{-*}*(qa+b) = q \Longleftrightarrow qa+b = q^2 c + q d \Longleftrightarrow b=c=0, d=a.$$
In particular, the faithful action of $GL(2,\hh)/N$ is not free. Moreover, the orbit-stabilizer theorem yields the following.

\begin{proposition}\label{representation}
In the set $\rft$ of regular fractional transformations, which is the orbit of the identity function in $\mathcal{L}_\infty$, two elements $\f_A, \f_B$ coincide if, and only if, there exists $c \in \hh \setminus \{0\}$ such that $B = cA$. In particular, for all $\f \in \rft$ either $\f$ is a regular affine transformation or there exist (unique) $a,b,p \in \hh$ such that  
\begin{equation}
\f(q) = (q-p)^{-*}*(qa + b).
\end{equation}
\end{proposition}

We conclude this section observing that proposition \ref{quotients} relates $f.A = (f c + d)^{-*}*(f a+ b)$ to $F_A \circ f = (f c + d)^{-1}(f a+ b)$ as follows.

\begin{remark}\label{relation}
Let $A =
\left( \begin{array}{cc}
a & c \\
b & d
\end{array}\right) \in GL(2,\hh)$ 
and let $f\in \mathcal{L}_R$. If $T = T_{f c + d}$ then
\begin{equation}
f.A = F_A \circ f \circ T
\end{equation}
wherever $(f c + d)^s$ does not vanish.
\end{remark}

\begin{proof}
Set $g(q) = f(q)c+d$ and $h(q) = f(q)a+b$ and notice that, by definition, $f.A = g^{-*}*h$. Now, according to proposition \ref{quotients}, $g^{-*}*h(q) = [g\circ T_g(q)]^{-1} [h \circ T_g(q)]$ with $T_g(q) = g^c(q)^{-1} q g^c(q)$. Now, $g\circ T_g(q) = f(T_g(q))c+d$ and $h \circ T_g(q) = f(T_g(q))a+b$ so that
$$f.A = g^{-*}*h = [(f\circ T_g)c+d]^{-1} [(f\circ T_g)a+b] = F_A \circ f \circ T_g$$
as desired.
\end{proof}


\section{Characterization of regular fractional transformations}\label{sectioncharacterization}

In section \ref{sectionregularfractional}, we treated regular fractional transformations as purely algebraic objects. We now study them as functions, determining their domains of definition and their properties.

The case of regular affine transformations  $\f(q) = qa + b \in \aff$ has already been studied: such an $\f$ is a biregular function $\hh \to \hh$ which extends to a homeomorphism from $\rh$ onto itself by setting $\f(\infty) = \infty$; we saw that $\f$ is also a linear fractional transformation.
Thanks to proposition \ref{representation}, all regular fractional transformations which are not affine are of the form $\f(q) = (q-p)^{-*}*(qa + b)$ for some $a,b,p \in \hh$. In order to study such transformations, we will make use of the following lemma, (see \cite{milan} for the proof).

\begin{lemma}\label{lemmatecnico} For any two quaternions $\alpha, \beta$ with $\beta
\neq \bar \alpha$ belonging to a same sphere $x+y\s$, we have $(\beta-\bar \alpha)^{-1}\beta(\beta-\bar \alpha)=\alpha$.
\end{lemma}

\begin{theorem}\label{diversificazionepoli} Consider the regular fractional transformation $\f(q) = (q-p)^{-*}*(qa + b)$ with $a,b,p \in \hh$ and let $F(q) = (q-p)^{-1}(q a + b)$. If $p=x\in \rr$ then $\f = F$ and in particular $\f$ extends to a homeomorphism from $\rh$ onto itself. If, on the contrary, $p=x+yI$ with $y\neq 0$ then $\f$ is a homeomorphism of $\rh  \setminus (x+y\s)$ onto $\rh \setminus F(x+y\s)$ and $\f$ does not extend to $\rh$ as a continuous function.
\end{theorem}
\begin{proof}The first statement follows by direct computation. As for the second one, $\f$ is related to the linear fractional transformation $F$ by remark \ref{relation}: $\f = F \circ T$, where  
$$T(q) = (q - \bar p)^{-1} q (q - \bar p)$$
is a homeomorphism from $\hh \setminus (x+y\s)$ onto itself. Remark that $T$ extends to a homeomorphism from $\rh \setminus (x+y\s)$ onto itself setting $T(\infty) = \infty$. By lemma \ref{lemmatecnico}, $T$ maps all points of $(x+y\s)\setminus\{\bar p\}$ to $p$. On the other hand the restriction of $T$ to the complex line $L_I$ through $p, \bar p$ coincides with the identity function of $L_I$, so that $\lim_{z \to \bar p, z \in L_I} T(z) = \bar p$. We conclude that $T$ does not extend to $\rh$ as a continuous function and, as a consequence, neither does $\f$.
\end{proof}

According to \cite{poincare}, $F(x+y\s)$ is either a 2-sphere or a plane. Since $F(x+yI) = \infty$, we conclude that $F(x+y\s)$ is a plane. Now let us prove the following characterization.

\begin{theorem}
Let $f : \rh \to \rh$ be a continuous, injective function. If $f$ is semiregular in $\hh$ then either $f$ is a regular affine transformation or it is a regular fractional transformation with $f(x) = \infty$ at a real point $x \in \rr$.
\end{theorem}

\begin{proof} If $f(\hh) \subseteq \hh$, then $f_{|_\hh} : \hh \to \hh$ is an injective regular map. By theorem \ref{affine}, $f$ must be a regular affine transformation.

If, on the contrary, $f$ equals $\infty$ at some point of $\hh$ then such point is unique by the injectivity of $f$. By theorem \ref{distribution}, this can only happen if the point is real, say $x \in \rr$. Setting $\tilde f (q) = f(q^{-1}+x)$ defines an injective function $\tilde f: \rh \to \rh$ which is semiregular in $\hh$. Since $\tilde f(\infty) = f(x) = \infty$ and $\tilde f$ is injective, $\tilde f$ cannot have a pole in $\hh = \rh \setminus \{\infty\}$. By the first part of the proof, $\tilde f$ is a regular affine transformation, say $\tilde f(q) =  q a + b$ for some $a,b \in \hh, a \neq 0$. Hence $f(q) = \tilde f((q-x)^{-1}) = (q-x)^{-1}a + b = (q-x)^{-1}(q b + a -xb) = (q-x)^{-*}*(q b + a -xb)$.
\end{proof}

The two theorems above yield:

\begin{corollary}\label{characterization}
Let $f$ be semiregular in $\hh$. Then $f$ extends to a homeomorphism from $\rh$ onto itself if and only if $f$ extends to a continuous, injective function in $\rh$, if and only if $f$ is a regular fractional transformation with pole in $\widehat{\rr}$, if and only if $f \in \rft \cap \lft$.
\end{corollary}

Remark that the set $\rft \cap \lft$ of regular fractional transformations which are also linear fractional is not a subgroup of $\lft$.


\section{Regular transformations of the unit ball}\label{sectionball}

We now focus on the transformations mapping the open unit ball $\B = B(0,1)$ onto itself. The set $\m = \{F \in \lft : F(\B) = \B\}$ of \emph{Moebius transformations} can be characterized as follows. Consider the symplectic group $Sp(1,1) = \{C \in GL(2, \hh) : \overline C ^t H C = H\}$ with $H = 
\left( \begin{array}{cr}
1 & 0 \\
0 & -1
\end{array}\right)$. One can prove, as in \cite{poincare}, that

\begin{theorem}
For all $A \in SL(2, \hh)$, the linear fractional transformation $F_A$ maps $\B$ onto itself if and only if $A \in Sp(1,1)$ if and only if there exist $u,v \in \partial \B, a \in \B$ such that
\begin{equation}
F_A(q) = v^{-1} (1-q \bar a)^{-1}(q-a)u
\end{equation}
for all $q \in \B$. In particular the antihomomorphism $\Phi$ defined in theorem \ref{antihomomorphism} can be restricted to a surjective group antihomomorphism $\Phi : Sp(1,1) \to \m$ with kernel $Ker(\Phi) = \{\pm \I\}$.
\end{theorem}

From this result we derive the following properties.

\begin{corollary}
For all $A \in SL(2, \hh)$, the regular fractional transformation $\f_A$ maps $\B$ onto itself if and only if $A \in Sp(1,1)$, if and only if there exist (unique) $u \in \partial \B, a \in \B$ such that
\begin{equation}
\f_A(q) = (1-q \bar a)^{-*}*(q-a)u
\end{equation}
for all $q \in \B$. In particular, the set $\mr = \{f \in \rft : f(\B) = \B\}$ of \emph{regular Moebius transformations} is the orbit of the identity function under the action of $Sp(1,1)$.
\end{corollary}

\begin{proof}
$\f_A$ maps $\B$ onto $\B$ if and only of $F_A$ does (indeed, by remark \ref{relation}, $\f_A = F_A \circ T$ where $T$ maps $\B$ onto itself). Moreover, $F_A(q) = v^{-1} (1-q \bar a)^{-1}(q-a)u = (v-q \bar a v)^{-1}(q-a)u$ implies 
$$\f_A(q) = (v-q \bar a v)^{-*}*(q-a)u =$$
$$= (1-q \bar v \bar a v)^{-*}*\bar v*(q-a)u = (1-q \bar v \bar a v)^{-*}*(q\bar v-\bar v a)u =$$
$$= (1-q \bar v \bar a v)^{-*}*(q-\bar v a v)\bar v u = (1-q \bar \alpha)^{-*}*(q-\alpha)\nu$$
with $\alpha = \bar v a v \in \B$ and $\nu = \bar v u \in \partial \B$.
\end{proof}

We will prove that all regular bijections $\B \to \B$ are regular Moebius transformations. We begin by rephrasing the Schwarz lemma (proven in \cite{advances}) as follows.

\begin{lemma}\label{schwarz}
Let $f : \B \to \B$ be a regular function such that $f(0) = 0$ and let $g : \B \to \hh$ be the regular function such that $f(q)= q*g(q) = q g(q)$. Then $|g(q)| \leq 1$ for all $q \in \B$ and equality holds at some point of $\B$ if and only if $g \equiv u$ for some $u \in \partial \B$.
\end{lemma}

Notice that the existence of $g$ is granted by theorem \ref{factorization}. We also need the following lemma.

\begin{lemma}\label{limit}
If $f : \B \to \B$ is regular and bijective then, for all $\nu \in \partial \B$, $\lim_{q \to \nu} |f(q)| = 1$.
\end{lemma}

\begin{proof}
Fix $\nu \in \partial \B$. Let $\{q_n\}_{n \in \nn}$ be a sequence converging to $\nu$, let $\mu$ be a limit point of $\{f(q_n)\}_{n \in \nn}$ and let us prove that $|\mu| = 1$. Suppose by contradiction that $|\mu|<1$. By the surjectivity of $f$, there exists $q \in \B$ such that $f(q) = \mu$. Choose an open neighborhood $U$ of $q$ in $\B$ such that $\nu \not \in \overline{U}$: there exists an $N \in \nn$ such that $q_n \not \in U$ for all $n \geq N$. By the injectivity of $f$, we conclude that $f(q_n) \not \in f(U)$ for $n\geq N$. But, according to the open mapping theorem \ref{open}, $f(U)$ is an open neighborhood of $\mu$. Thus $\mu$ cannot be a limit point for $\{f(q_n)\}_{n \in \nn}$, a contradiction with the hypothesis. 
\end{proof}

We are now ready to prove the desired result in one special case.

\begin{theorem}\label{rotation}
If $f : \B \to \B$ is regular and bijective and $f(0) = 0$ then there exists $u \in \partial \B$ such that $f (q) = q u$ for all $q \in \B$.
\end{theorem}

\begin{proof}
According to the Schwarz lemma \ref{schwarz}, there exists a regular function $g : \B \to \B$  such that $f(q)= q*g(q) = q g(q)$. Moreover, if we prove that $|g(0)| = 1$ then we can conclude that there exists $u \in \partial \B$ such that $g \equiv u$ as desired. 

Suppose by contradiction $|g(0)| < \varepsilon < 1$ for some $\varepsilon$. Since $\lim\limits_{q \to \partial \B} |g(q)|=1$ (thanks to lemma \ref{limit}), there exists a $\delta<1$ such that $|g(q)| \geq \varepsilon$ for $|q| \geq \delta$. The fact that $|g(0)| < \varepsilon \leq \min\limits_{|q| = \delta} |g(q)|$ implies that $g$ is not constant and has a minimum modulus point $p$ in $B(0,\delta)$. By the minimum modulus principle \ref{minimum}, $g(p) =0$. This implies $f(p) = 0 = f(0)$, a contradiction with the injectivity assumption for $f$.
\end{proof}

We now prove the main result of this section.

\begin{theorem}
If $f : \B \to \B$ is regular and bijective then there exist $a \in \B, u \in \partial \B$ such that $f (q) = (1-q \bar a)^{-*}*(q-a)u$ for all $q \in \B$.
\end{theorem}

\begin{proof}
Let $c = f(0) \in \B$, $C = 
\left( \begin{array}{cr}
1 & -\bar c \\
-c & 1
\end{array}\right)$ and $\tilde f = f.C = (1-f \bar c)^{-*}*(f-c)$.
The function $\tilde f$ is regular in $\B$ because $1-f \bar c$ does not have zeros in $\B$ (since $|c|<1, |f(q)|<1$ for $|q|<1$). Now recall that by remark \ref{relation}
$$\tilde f = f.C = F_C \circ f \circ T$$
where $T = T_{1-f\bar c}$ is a bijection of $\B$ onto $\B$. Since $f : \B \to \B$ is bijective and $F_C(q) = (1-q\bar c)^{-1}(q-c)$ is a Moebius transformation of $\B$, the function $\tilde f$ is a bijection of $\B$ onto $\B$, too. Moreover, 
$$\tilde f(0) =  F_C \circ f \circ T(0) =  F_C \circ f(0) = F_C(c) = 0.$$
By theorem \ref{rotation}, we conclude that $\tilde f(q) = q u$ for some $u \in \partial \B$. 
Now, $f = (f.C).C^{-1} = \tilde f. C^{-1}$ where $C^{-1}$ coincides, up to real multiplicative constant, with the matrix 
$\left( \begin{array}{cr}
1 & \bar c \\
c & 1
\end{array}\right)$. Thus 
$$f(q) = (1+\tilde f \bar c)^{-*}*(\tilde f + c) = (1+ q u \bar c)^{-*}*(q u + c) =$$
$$= (1+q \overline{c \bar u})^{-*}*(q + c \bar u) u,$$
as desired.
\end{proof}

Our final remark is the following.

\begin{proposition}
For all $a,b \in \B$ there exists a transformation $\f \in \mr$ mapping $a$ to $b$.
\end{proposition}

\begin{proof}
It suffices to set $\f = \f_{ABM}$ with $$A = 
\left( \begin{array}{cr}
1 & -\bar a \\
-a & 1
\end{array}\right), B = 
\left( \begin{array}{cr}
1 & |b| \\
|b| & 1
\end{array}\right), M = 
\left( \begin{array}{cr}
\frac{b}{|b|} & 0 \\
0 & 1
\end{array}\right).$$
Indeed, by direct computation $\f_A(q) = (1-q \bar a)^{-*}*(q-a)$ maps $a$ to $0$, $\f_{AB} = \left(1 + \f_A |b|\right)^{-*}*\left(\f_A + |b|\right) =  \left(1 + \f_A |b|\right)^{-1}\left(\f_A + |b|\right)$ maps $a$ to $|b|$ and $\f_{ABM} = \f_{AB} \cdot \frac{b}{|b|}$ maps $a$ to $b$.
\end{proof}


\bibliography{Moebius}

\bibliographystyle{abbrv}


\end{document}